\documentclass{article}
\usepackage{amssymb}

\usepackage{eurosym}
\usepackage{amsfonts}
\usepackage{amsmath}

\newtheorem{theorem}{Theorem}

\newtheorem{lemma}[theorem]{Lemma}

\newtheorem{remark}[theorem]{Remark}

\begin{document}

\title{Lyapunov-Sylvester Computational Method for Two-Dimensional
Boussinesq Equation}
\author{Abdelhamid Bezia$^{\ast }$,~Anouar Ben Mabrouk~and~Kamel Betina \\
$^{\ast }$Faculty of Mathematics, USTHB\\
E-mail: abdelhamid.bezia@gmail.com.}
\maketitle

\begin{abstract}
A numerical method is developed leading to algebraic systems based on
generalized Lyapunov-Sylvester operators to approximate the solution of
two-dimensional Boussinesq equation. It consists of an order reduction
method and a finite difference discretization. It is proved to be uniquely
solvable, stable and convergent by using Lyapunov criterion and manipulating
Lyapunov-Sylvester operators. Some numerical implementations are provided at
the end of the paper to validate the theoretical results.\newline
\end{abstract}

\section{Introduction}

In the present work we propose to serve of algebraic operators to
approximate the solutions of some PDEs such as Boussinesq one in higher
dimensions without adapting classical developments based on separation of
variables, radial solutions, ...etc. The crucial idea is to prove that even
simple methods of discretization of PDEs such as finite difference, finite
volumes, can be transformed to well adapted algebraic systems such as
Lyapunov-Sylvester ones developed hereafter. In the present paper,
fortunately, we were confronted with more complicated but fascinating method
to prove the invertibility of the algebraic operator yielded in the
numerical scheme. Instead of using classical methods such as tri-diagonal
transformations we applied a topological method to prove the invertibility.
This is good as it did not necessitate to compute eigenvalues and precisely
bounds/estimates of eigenvalues or direct inverses which remains a
complicated problem in general linear algebra and especially for generalized
Lyapunov-Sylvester operators. Recall that even though, bounds/estimates of
eigenvalues can already be efficient in studying stability. Recall also that
block tridiagonal systems for classical methods can be already used here
also and can be solved for example using iterative techniques, or highly
structured bandwidth solvers, Kronecker-product techniques, etc. These
methods have been subjects of more general discretizations. See \cite%
{ElMikkawy1}, \cite{ElMikkawy2}, \cite{ElMikkawy3}, \cite{El-Mikkawy-Atlan},
\cite{ElMikkawy-Karawia}, \cite{Jia-Li} for a review on tridiagonal and
block tridiagonal systems, their advantages as well as their disadvantages.
In the present paper our principal aim is to apply other algebraic methods
to investigates numerical solutions for PDEs in multi-dimensional spaces. We
aim to prove that Lyapunov-Syalvester operators can be good candidates for
such aim and that they may give best solvers compared to tridiagonal and/or
block tridiagonal ones. The present paper is devoted to the development of a
numerical method based on two-dimensional finite difference scheme to
approximate the solution of the nonlinear Boussinesq equation in $\mathbb{R}%
^{2}$ written on the form
\begin{equation}  \label{eqn1-1}
u_{tt}=\Delta
u+qu_{xxxx}+(u^{2})_{xx},\quad((x,y),t)\in\Omega\times(t_{0},+\infty)
\end{equation}%
with initial conditions
\begin{equation}  \label{eqn1-3}
u(x,y,t_{0})=u_{0}(x,y)\quad\hbox{and}\quad\displaystyle\frac{\partial\,u}{%
\partial\,t}(x,y,t_{0})=\varphi(x,y),\quad(x,y)\in\Omega
\end{equation}%
and boundary conditions
\begin{equation}  \label{eqn1-4}
\frac{\partial u}{\partial\eta}\left(x,y,t\right)=0,\quad((x,y),t)\in%
\partial\Omega\times(t_{0},+\infty).
\end{equation}
In order to reduce the derivation order, we set
\begin{equation}  \label{eqn1-2}
v=qu_{xx}+u^{2}.
\end{equation}%
We have to solve the system
\begin{equation}  \label{systemequation1}
\left\{
\begin{array}{l}
u_{tt}=\Delta u+v_{xx},\quad(x,y,t)\in\Omega\times(t_{0},+\infty) \\
v=qu_{xx}+u^{2},\quad(x,y,t)\in\Omega\times(t_{0},+\infty) \\
\left(u,v\right)(x,y,t_{0})=\left(u_{0},v_{0}\right)(x,y),\quad(x,y)\in%
\overline{\Omega} \\
\frac{\partial\,u}{\partial\,t}(x,y,t_{0})=\varphi(x,y),\quad(x,y)\in%
\overline{\Omega} \\
\frac{\partial }{\partial\eta}(u,v)(x,y,t)=0,\quad(x,y,t)\in\partial\Omega%
\times(t_{0},+\infty)%
\end{array}
\right.
\end{equation}%
on a rectangular domain $\Omega=]L_{0},L_{1}[\times]L_{0},L_{1}[$ in $%
\mathbb{R}^{2}$. $t_{0}\geq0$ is a real parameter fixed as the initial time,
$u_{tt}$ is the second order partial derivative in time, $\Delta=\frac{%
\partial^{2}}{\partial\,x^{2}}+\frac{\partial^{2}}{\partial\,y^{2}}$ is the
Laplace operator in $\mathbb{R}^{2}$, $q$ is a real constant, $u_{xx}$ and $%
u_{xxxx}$ are respectively the second order and the fourth order partial
derivative according to $x$. $\displaystyle\frac{\partial}{\partial\eta}$ is
the outward normal derivative operator along the boundary $\partial\Omega$.
Finally, $u$, $u_{0}$ and $\varphi$ are real valued functions with $u_{0}$
and $\varphi$ are $\mathcal{C}^{2}$ on $\overline{\Omega}$. $u$ (and
consequently $v$) is the unknown candidates supposed to be $\mathcal{C}^{4}$
on $\overline{\Omega}$.

The Boussinesq equation has a wide reputability in both theoretic and
applied fields such as hydrodynamics, traveling-waves. It governs the flow
of ground water, Heat conduction, natural convection in thermodynamics for
both volume and fluids in porous media, etc. For this reason, many studies
have been developed discussing the solvability of such equations. There are
works dealing with traveling-wave solutions, self-similar solutions,
scattering method, mono and multi dimensional versions, reduction of multi
dimensional equations with respect to algebras, etc. In \cite{Dehghan2006},
several finite difference schemes such as three fully implicit finite
difference schemes, two fully explicit finite difference techniques, an
alternating direction implicit procedure and the Barakat and Clark type
explicit formula are discussed and applied to solve the two-dimensional
Schrodinger equation with Dirichlets boundary conditions. In \cite%
{Dehghan2008}, the solution of a generalized Boussinesq equation has been
developed by means of the homotopy perturbation method. It consisted in a
technique method that avoids the discretization, linearization, or small
perturbations of the equation and thus reduces the numerical computations.
In \cite{Shokri-Dehghan}, a collocation and approximation of the numerical
solution of the improved Boussinesq equation is obtained based on radial
bases. A predic atorcorrector scheme is provided and the Not-a-Knot method
is used to improve the accuracy in the boundary. Next, \cite{Dehghan-Salehi}%
, a boundary-only meshfree method has been applied to approximate the
numerical solution of the classical Boussinesq equation in one dimension.
See for instance \cite%
{BaoDan1,Benmabrouk1,Bratsos1,Clarkson1,Jafarizadeh1,Kano1,Kaya1,Lai1,Liu1,Parlange1,Song1,Varlamov1,Wazwaz1,Yi1,Yi2}
and the references therein for backgrounds on theses facts.

The method developed in this paper consists in replacing time and space
partial derivatives by finite-difference approximations in order linear
Lyapunov systems. An order reduction method is adapted leading to a system
of coupled PDEs which is transformed by the next to a discrete algebraic
one. The motivation behind the idea of applying Lyapunov operators was
already evoked in our work \cite{Benmabrouk1}. We recall in brief that such
a method leads to fast convergent and more accurate discrete algebraic
systems without going back to the use of tridiagonal and/or
fringe-tridiagonal matrices already used when dealing with multidimensional
problems especially in discrete PDEs.

In the organization of the present paper, the next section is concerned with
the introduction of the finite difference scheme. Section 3 is devoted to
the discretization of the continuous reduced system obtained from (\ref%
{eqn1-1})-(\ref{eqn1-4}) by the order reduction method. Section 4 deals with
the solvability of the discrete Lyapunov equation obtained from the
discretization developed in the section 3. In section 5, the consistency of
the method is shown and next, the stability and convergence of are proved
based on Lyapunov method. Finally, a numerical implementation is provided in
section 6 leading to the computation of the numerical solution and error
estimates.

\section{Discrete Two-Dimensional Boussinesq Equation}

Consider the domain $\Omega =]L_{0},L_{1}[\times]L_{0},L_{1}[\subset \mathbb{%
R}^{2}$ and an integer $J\in\mathbb{N}^{\ast}$. Denote $h=\displaystyle\frac{%
L_{1}-L_{0}}{J}$ for the space step, $x_{j}=L_{0}+jh$ and $y_{m}=L_{0}+mh$
for all $(j,m)\in {I}^{2}=\{0,1,\dots,J\}^{2}$. Let $l=\Delta\,t$ be the
time step and $t_{n}=t_{0}+nl$, $n\in\mathbb{N}$ for the discrete time grid.
For $(j,m)\in {I}$ and $n\geq 0$, $u_{j,m}^{n}$ will be the net function $%
u(x_{j},y_{m},t_{n})$ and $U_{j,m}^{n} $ the numerical solution. The
following discrete approximations will be applied for the different
differential operators involved in the problem. For time derivatives, we set
\begin{equation*}
\displaystyle\,u_{t}\rightsquigarrow\displaystyle\frac{%
U_{j,m}^{n+1}-U_{j,m}^{n-1}}{2l}\quad\hbox{and}\quad\displaystyle%
\,u_{tt}\rightsquigarrow \displaystyle\frac{%
U_{j,m}^{n+1}-2U_{j,m}^{n}+U_{j,m}^{n-1}}{l^{2}}
\end{equation*}%
and for space derivatives, we shall use
\begin{equation*}
\displaystyle\,u_{x}\rightsquigarrow\displaystyle\frac{%
U_{j+1,m}^{n}-U_{j-1,m}^{n}}{2h}\quad \hbox{and}\quad\displaystyle %
\,u_{y}\rightsquigarrow \displaystyle\frac{U_{j,m+1}^{n}-U_{j,m-1}^{n}}{2h}
\end{equation*}%
for first order derivatives and
\begin{equation*}
\displaystyle\,u_{xx}\rightsquigarrow\displaystyle\frac{U_{j+1,m}^{n,%
\alpha}-2U_{j,m}^{n,\alpha}+U_{j-1,m}^{n,\alpha}}{h^{2}},\quad\displaystyle %
\,u_{yy}\rightsquigarrow\displaystyle\frac{U_{j,m+1}^{n,\alpha}-2U_{j,m}^{n,%
\alpha}+U_{j,m-1}^{n,\alpha}}{h^{2}}
\end{equation*}%
for second order ones, where for $n\in\mathbb{N}^{\ast}$ and $\alpha\in%
\mathbb{R}$,
\begin{equation*}
u^{n,\alpha}=\alpha\,U^{n+1}+(1-2\alpha)U^{n}+\alpha\,U^{n-1}.
\end{equation*}%
Finally, we denote $\sigma=\displaystyle\frac{l^{2}}{h^{2}}$ and $\delta=%
\displaystyle\frac{q}{h^{2}}$.

For $(j,m)\in \mathring{I}^{2}$ an interior point of the grid $I^{2}$, ($%
\mathring{I}=\{1,2,\dots ,J-1\}$), and $n\geq 1$, the following discrete
equation is deduced from the first equation in the system (\ref%
{systemequation1}).
\begin{equation}
\begin{matrix}
&  & U_{j,m}^{n+1}-2U_{j,m}^{n}+U_{j,m}^{n-1}\hfill\cr & = &
\sigma\alpha%
\left(U_{j-1,m}^{n+1}-4U_{j,m}^{n+1}+U_{j+1,m}^{n+1}+U_{j,m-1}^{n+1}+U_{j,m+1}^{n+1}\right)\hfill%
\cr &  & +\sigma(1-2\alpha)%
\left(U_{j-1,m}^{n}-4U_{j,m}^{n}+U_{j+1,m}^{n}+U_{j,m-1}^{n}+U_{j,m+1}^{n}%
\right)\hfill\cr &  & +\sigma\alpha%
\left(U_{j-1,m}^{n-1}-4U_{j,m}^{n-1}+U_{j+1,m}^{n-1}+U_{j,m-1}^{n-1}+U_{j,m+1}^{n-1}\right)\hfill%
\cr &  & +\sigma\alpha\left(V_{j-1,m}^{n+1}-2V_{j,m}^{n+1}+V_{j+1,m}^{n+1}%
\right)\hfill\cr &  & +\sigma(1-2\alpha)%
\left(V_{j-1,m}^{n}-2V_{j,m}^{n}+V_{j+1,m}^{n}\right)\hfill\cr &  &
+\sigma\alpha\left(V_{j-1,m}^{n-1}-2V_{j,m}^{n-1}+V_{j+1,m}^{n-1}\right).%
\hfill%
\end{matrix}
\label{eqn1-1discrete}
\end{equation}
Similarly, the following discrete equation is obtained from equation (\ref%
{eqn1-2}).
\begin{equation}
\begin{matrix}
\,V_{j,m}^{n+1}+V_{j,m}^{n-1} & = & 2\delta \alpha \left(
U_{j-1,m}^{n+1}-2U_{j,m}^{n+1}+U_{j+1,m}^{n+1}\right) \hfill \cr &  &
+2\delta (1-2\alpha )\left( U_{j-1,m}^{n}-2U_{j,m}^{n}+U_{j+1,m}^{n}\right)
\hfill \cr &  & +2\delta \alpha \left(
U_{j-1,m}^{n-1}-2U_{j,m}^{n-1}+U_{j+1,m}^{n-1}\right) \hfill \cr &  & +2%
\widehat{F(U_{j,m}^{n})}\hfill%
\end{matrix}
\label{eqn2-1discrete}
\end{equation}%
where
\begin{equation*}
F(u)=u^{2},\;F^{n}=F(u^{n})\;\;\hbox{and}\;\;\widehat{F^{n}}=\displaystyle%
\frac{F^{n-1}+F^{n}}{2}.
\end{equation*}
The discrete boundary conditions are written for $n\geq 0 $ as
\begin{equation}
\displaystyle\,U_{1,m}^{n}=U_{-1,m}^{n}\quad \hbox{and}\quad \displaystyle%
\,U_{J-1,m}^{n}=U_{J+1,m}^{n},  \label{CBD1}
\end{equation}%
\begin{equation}
\displaystyle\,U_{j,1}^{n}=U_{j,-1}^{n}\quad \hbox{and}\quad \displaystyle%
\,U_{j,J-1}^{n}=U_{j,J+1}^{n},  \label{CBD2}
\end{equation}
The parameter $q$ is related to the equation and has the role of a
viscosity-type coefficient and thus it is related to the physical domain of
the model. The barycenter parameter $\alpha$ is used to calibrates the
position of the approximated solution around the exact one. Of course, these
parameters affect surely the numerical solution as well as the error
estimates. This fact will be recalled later in the numerical implementations
part. In a future work in progress now, we are developing results on
numerical solutions of 2D Schr\"odinger equation on the error estimates as a
function on the barycenter calibrations by using variable coefficients $%
\alpha_n$ instead of constant $\alpha$. The use of these calibrations
permits the use of implicit/explicit schemes by using suitable values. For
example for $\alpha=\frac{1}{2}$, the barycentre estimation
\begin{equation*}
V^{n,\alpha}=\alpha\,V^{n+1}+(1-2\alpha)V^{n}+\alpha\,V^{n-1}=\displaystyle%
\frac{V^{n+1}+V^{n-1}}{2}
\end{equation*}
which is an implicit estimation that guarantees an error of order $2$ in
time.

Next, as it is motioned in the introduction, the idea consists in applying
Lyapunov-Sylvester operators to approximate the solution of the continuous
problem (\ref{eqn1-1})-(\ref{eqn1-4}) or its discrete equivalent system (\ref%
{eqn1-1discrete})-(\ref{CBD2}). Denote
\begin{equation*}
a_{1}=\displaystyle\frac{1}{2}+2\alpha\sigma,\quad\,a_{2}=-\alpha\sigma,
\end{equation*}%
\begin{equation*}
b_{1}=1-2(1-2\alpha)\sigma,\quad\,b_{2}=(1-2\alpha)\sigma,
\end{equation*}%
\begin{equation*}
c_{1}=(1-2\alpha)\delta\quad\hbox{and}\quad\,c_{2}=\alpha\delta.
\end{equation*}%
Equation (\ref{eqn1-1discrete}) becomes
\begin{equation}  \label{eqn1-1discrete-forme-a}
\begin{matrix}
&  &
a_{2}U_{j-1,m}^{n+1}+a_{1}U_{j,m}^{n+1}+a_{2}U_{j+1,m}^{n+1}+a_{2}U_{j,m-1}^{n+1}+a_{1}U_{j,m}^{n+1}+a_{2}U_{j,m+1}^{n+1}\hfill%
\cr &  & +a_{2}\left(V_{j-1,m}^{n+1}-2V_{j,m}^{n+1}+V_{j+1,m}^{n+1}\right)%
\hfill\cr & = &
b_{2}U_{j-1,m}^{n}+b_{1}U_{j,m}^{n}+b_{2}U_{j+1,m}^{n}+b_{2}U_{j,m-1}^{n}+b_{1}U_{j,m}^{n}+b_{2}U_{j,m+1}^{n}\hfill%
\cr &  &
-a_{2}U_{j-1,m}^{n-1}-a_{1}U_{j,m}^{n-1}-a_{2}U_{j+1,m}^{n-1}-a_{2}U_{j,m-1}^{n-1}-a_{1}U_{j,m}^{n-1}-a_{2}U_{j,m+1}^{n-1}\hfill%
\cr &  & +b_{2}\left(V_{j-1,m}^{n}-2V_{j,m}^{n}+V_{j+1,m}^{n}\right)\hfill\cr
&  & -a_{2}\left(V_{j-1,m}^{n-1}-2V_{j,m}^{n-1}+V_{j+1,m}^{n-1}\right).\hfill%
\end{matrix}%
\end{equation}
Equation (\ref{eqn2-1discrete}) becomes
\begin{equation}  \label{eqn2-1discrete-forme-a}
\begin{matrix}
&  & V_{j,m}^{n+1}-2c_{2}%
\left(U_{j-1,m}^{n+1}-2U_{j,m}^{n+1}+U_{j+1,m}^{n+1}\right)\hfill\cr & = &
2c_{1}\left(U_{j-1,m}^{n}-2U_{j,m}^{n}+U_{j+1,m}^{n}\right)\hfill\cr &  &
+2c_{2}\left(U_{j-1,m}^{n-1}-2U_{j,m}^{n-1}+U_{j+1,m}^{n-1}\right)\hfill\cr
&  & -V_{j,m}^{n-1}+2\widehat{F(U_{j,m}^{n})}.\hfill%
\end{matrix}%
\end{equation}
Denote $A$, $B$ and $R$ the matrices defined by
\begin{equation*}
A=%
\begin{pmatrix}
a_{1} & 2a_{2} & 0 & ... & ... & 0\hfill\cr a_{2} & a_{1} & a_{2} & \ddots &
\ddots & \vdots\hfill\cr 0 & \ddots & \ddots & \ddots & \ddots & \vdots\hfill%
\cr \vdots & \ddots & \ddots & \ddots & \ddots & 0\hfill\cr \vdots & \ddots
& \ddots & a_{2} & a_{1} & a_{2}\hfill\cr 0 & ... & ... & 0 & 2a_{2} &
a_{1}\hfill%
\end{pmatrix}
,\quad B=
\begin{pmatrix}
b_{1} & 2b_{2} & 0 & ... & ... & 0\hfill\cr b_{2} & b_{1} & b_{2} & \ddots &
\ddots & \vdots\hfill\cr 0 & \ddots & \ddots & \ddots & \ddots & \vdots\hfill%
\cr \vdots & \ddots & \ddots & \ddots & \ddots & 0\hfill\cr \vdots & \ddots
& \ddots & b_{2} & b_{1} & b_{2}\hfill\cr 0 & ... & ... & 0 & 2b_{2} &
b_{1}\hfill%
\end{pmatrix}%
\end{equation*}
and
\begin{equation*}
R=
\begin{pmatrix}
-2 & 2 & 0 & ... & ... & 0\hfill \cr1 & -2 & 1 & \ddots & \ddots & \vdots
\hfill \cr0 & \ddots & \ddots & \ddots & \ddots & \vdots \hfill \cr\vdots &
\ddots & \ddots & \ddots & \ddots & 0\hfill \cr\vdots & \ddots & \ddots & 1
& -2 & 1\hfill \cr0 & ... & ... & 0 & 2 & -2\hfill%
\end{pmatrix}%
\end{equation*}
The system (\ref{CBD1})-(\ref{eqn2-1discrete-forme-a}) can be written on the
matrix form
\begin{equation}
\left\{
\begin{matrix}
\mathcal{L}_{A}(U^{n+1})+a_{2}RV^{n+1}=\mathcal{L}_{B}(U^{n})-\mathcal{L}%
_{A}(U^{n-1})+R(b_{2}V^{n}-a_{2}V^{n-1}),\hfill\cr V^{n+1}-2c_{2}RU^{n+1}\!=%
\!2R(c_{1}U^{n}+c_{2}U^{n-1})-V^{n-1}\!+\!2\widehat{F^{n}}\hfill%
\end{matrix}
\right.  \label{matrixform1}
\end{equation}%
for all $n\geq1$ where
\begin{equation*}
U^{n}=(U_{j,m}^{n})_{0\leq\,j,m\leq\,J},\;V^{n}=(V_{j,m}^{n})_{0\leq\,j,m%
\leq\,J}\;\hbox{and}\;F^{n}=(F(U_{j,m}^{n}))_{0\leq\,j,m\leq\,J}
\end{equation*}
and for a matrix $Q\in\mathcal{M}_{(J+1)^{2}}(\mathbb{R})$, $\mathcal{L}_{Q}$
is the Lyapunov operator defined by
\begin{equation*}
\mathcal{L}_{Q}(X)=QX+XQ^{T},\;\forall\,X\in\mathcal{M}_{(J+1)^{2}}(\mathbb{R%
}).
\end{equation*}
Remark that $V$ is obtained from the auxiliary function $v$ that is applied
to reduce the order of the original PDEs in $u$. This reduction yielded the
Lyapunov-Syslvester system (\ref{matrixform1}) above. A natural question
that can be raised here turns around the ordering of $U$ and $V$. So, we
stress the fact that no essential idea is fixed at advance but, this is
strongly related to the system obtained. For example, in (\ref{matrixform1})
above, it is easy to substitute the second equation into the first to omit
the unknown matrix $V^{n+1}$ from the first equation. But in the contrary,
it is not easier to do the same for $U^{n+1}$, due to the difficulty to
substitute it from $\mathcal{L}_{A}(U^{n+1})$. It is also not guaranteed
that the part $a_{2}RV^{n+1}$ in the first equation is invertible to
substitute $V^{n+1}$. So, it is essentially the final system that shows the
ordering in $U$ and $V$.

\section{Solvability of the Discrete Problem}

In \cite{Benmabrouk1}, the authors have transformed the Lyapunov operator
obtained from the discretization method into a standard linear operator
acting on one column vector by juxtaposing the columns of the matrix $X$
horizontally which leads to an equivalent linear operator characterized by a
fringe-tridiagonal matrix. We used standard computation to prove the
invertibility of such an operator. Here. we do not apply the same
computations as in \cite{Benmabrouk1}, but we develop different arguments.
The first main result is stated as follows.

\begin{theorem}
\label{theorem1} The system (\ref{matrixform1}) is uniquely solvable
whenever $U^{0}$ and $U^{1}$ are known.
\end{theorem}

\textit{Proof.} It reposes on the inverse of Lyapunov operators. Consider
the endomorphism $\Phi $ defined on $\mathcal{M}_{(J+1)^{2}}(\mathbb{R}%
)\times \mathcal{M}_{(J+1)^{2}}(\mathbb{R}) $ by $\Phi
(X,Y)=(AX+XA^T+a_{2}RY,\displaystyle\frac{1}{2}Y-c_{2}RX)$. To prove Theorem %
\ref{theorem1}, it suffices to show that $ker\Phi $ is reduced to $0$.
Indeed,
\begin{equation*}
\Phi (X,Y)=0\Longleftrightarrow (AX+XA^T+a_{2}RY,\displaystyle\frac{1}{2}%
Y-c_{2}RX)=(0,0)
\end{equation*}%
or equivalently,
\begin{equation*}
Y=2c_{2}RX\quad \hbox{and}\quad \,(A+2a_{2}c_{2}R^{2})X+XA^T=0.
\end{equation*}%
So, the problem is transformed to the resolution of a Lyapunov type equation
of the form
\begin{equation}
\mathcal{L}_{W,A}(X)=WX+XA^{T}=0  \label{Lyapunovequation}
\end{equation}%
where $W$ is the matrix given by $W=A+2a_{2}c_{2}R^{2}$. Denoting
\begin{equation*}
\omega=2a_{2}c_{2},\,\,\omega_{1}=a_{1}+6\omega,\;\;\overline{\omega}%
_{1}=\omega_{1}+\omega\quad\hbox{and}\quad\omega_{2}=a_{2}-4\omega
\end{equation*}%
the matrix $W$ is explicitly given by
\begin{equation*}
W=
\begin{pmatrix}
\omega _{1} & 2\omega _{2} & 2\omega & 0 & \dots & \dots & \dots & 0 &
\hfill \cr\omega _{2} & \overline{\omega }_{1} & \omega _{2} & \omega &
\ddots & \ddots & \ddots & \vdots & \hfill \cr\omega & \omega _{2} & \omega
_{1} & \omega _{2} & \omega & \ddots & \ddots & \vdots & \hfill \cr0 & \ddots
& \ddots & \ddots & \ddots & \ddots & \ddots & \vdots & \hfill \cr\vdots &
\ddots & \ddots & \ddots & \ddots & \ddots & \ddots & 0 & \hfill \cr\vdots &
\ddots & \ddots & \omega & \omega _{2} & \omega _{1} & \omega _{2} & \omega
& \hfill \cr\vdots & \ddots & \ddots & \ddots & \omega & \omega _{2} &
\overline{\omega }_{1} & \omega _{2} & \hfill \cr0 & \dots & \dots & \dots &
0 & 2\omega & 2\omega _{2} & \omega _{1} & \hfill%
\end{pmatrix}%
\end{equation*}%
Next, we use the following preliminary result of differential calculus (See
\cite{HenriCartan} for example).

\begin{lemma}
\label{LemmeInversion} Let $E$ be a finite dimensional ($\mathbb{R}$ or $%
\mathbb{C}$) vector space and $(\Phi_{n})_{n}$ be a sequence of
endomorphisms converging uniformly to an invertible endomorphism $\Phi$.
Then, there exists $n_{0}$ such that, for any $n\geq\,n_{0}$, the
endomorphism $\Phi_{n}$ is invertible.
\end{lemma}

The proof is simple and can be found anywhere in differential calculus
references such as \cite{HenriCartan}. We recall it here for the convenience
and clearness of the paper. Recall that the set $Isom(E)$ (the set of
isomorphisms on $E$) is already open in $L(E)$ (the set of endomorphisms of $%
E$). Hence, as $\Phi\in\,Isom(E)$ there exists a ball $B(\Phi,r)\subset%
\,Isom(E)$. The elements $\Phi_n$ are in this ball for large values of $n$.
So these are invertible.

Assume now that $l=o(h^{2+s})$, with $s>0$ which is always possible. Then,
the coefficients appearing in $A$ and $W$ will satisfy as $h\longrightarrow0$
the following.
\begin{equation*}
A_{i,i}=\displaystyle\frac{1}{2}+\varepsilon\,h^{2+2s}\longrightarrow%
\displaystyle\frac{1}{2}.
\end{equation*}
For $1\leq\,i\leq\,J-1$,
\begin{equation*}
A_{i,i-1}=A_{i,i+1}=\displaystyle\frac{A_{0,1}}{2}=\displaystyle\frac{%
A_{J,J-1}}{2}=-\varepsilon\,h^{2+2s}\longrightarrow0.
\end{equation*}
For $2\leq\,i\leq\,J-2$,
\begin{equation*}
W_{i,i}=W_{0,0}=W_{J,J}=\displaystyle\frac{1}{2}+2\alpha\varepsilon%
\,h^{2+2s}-12\alpha^2\varepsilon\,h^{2s}\longrightarrow\displaystyle\frac{1}{%
2}.
\end{equation*}
Similarly,
\begin{equation*}
W_{1,1}=W_{J-1,J-1}=\displaystyle\frac{1}{2}+2\alpha\varepsilon\,h^{2+2s}-14%
\alpha^2\varepsilon\,h^{2s}\longrightarrow\displaystyle\frac{1}{2}
\end{equation*}
and
\begin{equation*}
W_{i,i-1}=W_{i,i+1}=\displaystyle\frac{W_{0,1}}{2}=\displaystyle\frac{%
W_{J,J-1}}{2}=-\alpha\varepsilon\,h^{2+2s}+8\alpha^2\varepsilon\,h^{2s}%
\longrightarrow0
\end{equation*}
Finally,
\begin{equation*}
W_{i,i-2}=W_{i,i+2}=\displaystyle\frac{W_{0,2}}{2}=\displaystyle\frac{%
W_{J,J-2}}{2}=-2\alpha^2\varepsilon\,h^{2s}\longrightarrow0.
\end{equation*}
Recall that the technique assumption $l=o(h^{2+s})$ is a necessary
requirement for the resolution of the present problem and may not be
necessary in other PDEs. See for example \cite{Benmabrouk2}, \cite%
{Benmabrouk1} and \cite{Benmabrouk3} for NLS and Heat equations. Next,
observing that for all $X$ in the space $\mathcal{M}_{(J+1)^2}(\mathbb{R}%
)\times\mathcal{M}_{(J+1)^2}(\mathbb{R})$,
\begin{equation*}
\begin{matrix}
\|(\mathcal{L}_{W,A}-I)(X)\| & =\|(W-\frac{1}{2}I)X+X(A^T-\frac{1}{2}%
I)\|\hfill \cr & \leq\left[\|W-\frac{1}{2}I\|+\|A^T-\frac{1}{2}I\|\right]%
\|X\|,\hfill%
\end{matrix}%
\end{equation*}
it results that
\begin{equation}  \label{LWAtendsVersId}
\|\mathcal{L}_{W,A}-I\|\leq\|W-\displaystyle\frac{1}{2}I\|+\|A^T-%
\displaystyle\frac{1}{2}I\|\leq\,C(\alpha)h^{2s}.
\end{equation}
Consequently, the Lyapunov endomorphism $\mathcal{L}_{W,A}$ converges
uniformly to the identity $I$ as $h$ goes towards 0 and $l=o(h^{2+s})$ with $%
s>0$. Using Lemma \ref{LemmeInversion}, the operator $\mathcal{L}_{W,A}$ is
invertible for $h$ small enough.

\section{Consistency, Stability and Convergence of the Discrete Method}

The consistency of the proposed method is done by evaluating the local
truncation error arising from the discretization of the system
\begin{equation}  \label{system1}
\left\{
\begin{matrix}
u_{tt}-\Delta\,u-v_{xx}=0,\hfill\cr v=qu_{xx}+u^2.\hfill%
\end{matrix}
\right.
\end{equation}
The principal part of the first equation is
\begin{equation}  \label{consistency1}
\begin{matrix}
\mathcal{L}_{u,v}^1(t,x,y) & = & \displaystyle\frac{l^2}{12}\displaystyle
\frac{\partial^4u}{\partial\,t^4} -\displaystyle\frac{h^2}{12}\left( %
\displaystyle\frac{\partial^4u}{\partial\,x^4}+\displaystyle\frac{%
\partial^4u }{\partial\,y^4}\right) -\alpha\,l^2\displaystyle\frac{%
\partial^2(\Delta\,u) }{\partial\,t^2}\hfill\cr & - & \displaystyle\frac{h^2%
}{12}\displaystyle \frac{\partial^2v}{\partial\,x^4} -\alpha\,l^2%
\displaystyle\frac{\partial^4v }{\partial\,t^2\partial\,x^2}%
+O(l^2+h^2).\hfill%
\end{matrix}%
\end{equation}
The principal part of the local error truncation due to the second part is
\begin{equation}  \label{consistency2}
\begin{matrix}
\mathcal{L}_{u,v}^2(t,x,y) & = & \displaystyle\frac{l^2}{2}\displaystyle
\frac{\partial^2v}{\partial\,t^2} +\displaystyle\frac{l^4}{24}\displaystyle
\frac{\partial^4v}{\partial\,t^4} -\displaystyle\frac{h^2}{12}\displaystyle
\frac{\partial^4u}{\partial\,x^4}\hfill\cr & - & \alpha\,l^2\displaystyle
\frac{\partial^4u}{\partial\,t^2\partial\,x^2}+O(l^2+h^2).\hfill%
\end{matrix}%
\end{equation}
It is clear that the two operators $\mathcal{L}_{u,v}^1$ and $\mathcal{L}%
_{u,v}^2$ tend toward 0 as $l$ and $h$ tend to 0, which ensures the
consistency of the method. Furthermore, the method is consistent with an
order 2 in time and space.

We now proceed by proving the stability of the method by applying the
Lyapunov criterion. A linear system $\mathcal{L}(x_{n+1},x_{n},x_{n-1},%
\dots)=0$ is stable in the sense of Lyapunov if for any bounded initial
solution $x_{0}$ the solution $x_{n}$ remains bounded for all $n\geq0$.
Here, we will precisely prove the following result.

\begin{lemma}
\label{LyapunovStabilityLemma} $\mathcal{P}_{n}$: The solution $%
(U^{n},V^{n}) $ is bounded independently of $n$ whenever the initial
solution $(U^{0},V^{0})$ is bounded.
\end{lemma}

We will proceed by recurrence on $n$. Assume firstly that $%
\|(U^0,V^0)\|\leq\eta$ for some $\eta$ positive. Using the system (\ref%
{matrixform1}), we obtain
\begin{equation}  \label{LyapunovStability2}
\left\{
\begin{matrix}
\mathcal{L}_{W,A}(U^{n+1})=\mathcal{L}_{\widetilde{B},B}(U^{n})+b_2RV^n-
\mathcal{L}_{W,A}(U^{n-1})-a_2R(F^{n-1}+F^n),\hfill\cr %
V^{n+1}=2c_2RU^{n+1}+2R(c_1U^{n}+c_2U^{n-1})-V^{n-1}+2\widehat{F^n}.\hfill%
\end{matrix}
\right.
\end{equation}
where $\widetilde{B}=B-2a_2c_1R^2$. Consequently,
\begin{equation}  \label{LyapunovStability3}
\begin{matrix}
\|\mathcal{L}_{W,A}(U^{n+1})\|\leq\|\mathcal{L}_{\widetilde{B}%
,B}\|.\|U^{n}\|+2|b_2|.\|V^n\|\hfill\cr \qquad\qquad\qquad\qquad\qquad +\|%
\mathcal{L}_{W,A}\|.\|U^{n-1}\|+2|a_2|(\|F^{n-1}\|+\|F^n\|)\hfill%
\end{matrix}%
\end{equation}
and
\begin{equation}  \label{LyapunovStability4}
\begin{matrix}
\|V^{n+1}\|\leq4|c_2|.\|U^{n+1}\|+4(|c_1|.\|U^{n}\|+|c_2|.\|U^{n-1}\|)\hfill%
\cr \qquad\qquad\qquad\qquad\qquad +\|V^{n-1}\|+\|F^{n-1}\|+\|F^n\|.\hfill%
\end{matrix}%
\end{equation}
Next, recall that, for $l=o(h^{s+2})$ small enough, $s>0$, we have
\begin{equation*}
a_1=\displaystyle\frac{1}{2}+2\alpha\,h^{2s+2}\rightarrow\displaystyle\frac{1%
}{2},\quad\,a_2=-\alpha\,h^{2s+2}\rightarrow0,
\end{equation*}
\begin{equation*}
b_1=1-2(1-2\alpha)h^{2s+2}\rightarrow1,\quad\,b_2=(1-2\alpha)h^{2s+2}%
\rightarrow0,
\end{equation*}
\begin{equation*}
c_1=(1-2\alpha)h^{-2}\rightarrow\infty\quad\hbox{and}\quad\,c_2=\alpha%
\,h^{2s+2}\rightarrow\infty,
\end{equation*}
\begin{equation*}
a_2c_1=-\alpha(1-2\alpha)h^{2s}\rightarrow0.
\end{equation*}
As a consequence, for $h$ small enough,
\begin{equation}  \label{LtildeBB}
\|\mathcal{L}_{\widetilde{B},B}\|\leq2\|B\|+2|a_2c_1|\|R\|^2\leq2%
\max(|b_1|,2|b_2|)+4|a_2c_1|\leq2+4=6,
\end{equation}
and the following lemma deduced from (\ref{LWAtendsVersId}).

\begin{lemma}
\label{LWABounded} For $h$ small enough, it holds for all $X\in \mathcal{M}%
_{(J+1)^{2}}(\mathbb{R})$ that
\begin{equation*}
\displaystyle\frac{1}{2}\|X\|\leq(1-C(\alpha)h^{2s})\|X\|\leq\|\mathcal{L}%
_{W,A}(X)\|\leq(1+C(\alpha)h^{2s})\|X\|\leq\displaystyle\frac{3}{2}\|X\|.
\end{equation*}
\end{lemma}

\hskip-17pt Indeed, recall that equation (\ref{LWAtendsVersId}) affirms that
$\|\mathcal{L}_{W,A}-I\|\leq\,C(\alpha)h^{2s}$ for some constant $C(\alpha)>0
$. Consequently, for any $X$ we get
\begin{equation*}
(1-C(\alpha)h^{2s})\|X\|\leq\|\mathcal{L}_{W,A}(X)\|\leq(1+C(\alpha)h^{2s})%
\|X\|.
\end{equation*}
For $h\leq\displaystyle\frac{1}{(C(\alpha))^{1/2s}}$, we obtain
\begin{equation*}
\displaystyle\frac{1}{2}\leq(1-C(\alpha)h^{2s})<(1+C(\alpha)h^{2s})\leq%
\displaystyle\frac{3}{2}
\end{equation*}
and thus Lemma \ref{LWABounded}. As a result, (\ref{LyapunovStability3})
yields that
\begin{equation}  \label{LyapunovStability3-1}
\displaystyle\frac{1}{2}\|U^{n+1}\|\leq6\|U^{n}\|+2\|V^{n}\|+\displaystyle%
\frac{3}{2}\|U^{n-1}\|+2(\|F^{n-1}\|+\|F^{n}\|).
\end{equation}%
For $n=0$, this implies that
\begin{equation}
\|U^{1}\|\leq12\|U^{0}\|+4\|V^{0}\|+3\|U^{-1}\|+4(\|F^{-1}\|+\|F^{0}\|).
\label{LyapunovStability3-2}
\end{equation}%
Using the discrete initial condition
\begin{equation*}
U^{0}=U^{-1}+l\varphi .
\end{equation*}%
Here we identify the function $\varphi $ to the matrix whom coefficients are
$\varphi_{j,m}=\varphi(x_{j},y_{m})$. We obtain
\begin{equation}
\|U^{-1}\|\leq\|U^{0}\|+l\|\varphi\|.  \label{U-1Bounds}
\end{equation}%
Observing that
\begin{equation*}
F_{j,m}^{-1}=F(U_{j,m}^{-1})=(U_{j,m}^{0}-l\varphi_{j,m})^{2},
\end{equation*}%
it results that
\begin{equation*}
|F_{j,m}^{-1}|\leq|U_{j,m}^{0}|^{2}+2l|\varphi_{j,m}|.|U_{j,m}^{0}|+l^{2}|%
\varphi_{j,m}|^{2}
\end{equation*}%
and consequently,
\begin{equation}
\|F^{-1}\|\leq\|U^{0}\|^{2}+2l\|\varphi\|.\|U^{0}\|+l^{2}\|\varphi\|^{2}.
\label{F-1Bounds}
\end{equation}%
Hence, equation (\ref{LyapunovStability3-2}) yields that
\begin{equation}
\|U^{1}\|\leq(15+8l\|\varphi\|)\|U^{0}\|+4\|V^{0}\|+8\|F^{0}\|+3l\|\varphi%
\|+4l^{2}\|\varphi\|^{2}.  \label{LyapunovStability3-3}
\end{equation}%
Now, the Lyapunov criterion for stability states exactly that
\begin{equation}  \label{LyapunovStability1}
\forall\,\,\varepsilon>0,\,\exists\,\eta>0\,\,\,s.t;\,\,\|(U^{0},V^{0})\|%
\leq\eta\,\,\Rightarrow\,\,\|(U^{n},V^{n})\|\leq\varepsilon,\,\,\forall\,n%
\geq0.
\end{equation}
For $n=1$ and $\|(U^{1},V^{1})\|\leq\varepsilon$, we seek an $\eta>0$ for
which $\|(U^{0},V^{0})\|\leq\eta$. Indeed, using (\ref{LyapunovStability3-3}%
), this means that, it suffices to find $\eta $ such that
\begin{equation}
8\eta^{2}+(19+8l\|\varphi\|)\eta+3l\|\varphi\|+4l^{2}\|\varphi\|^{2}-%
\varepsilon<0.  \label{LyapunovStability3-4}
\end{equation}%
The discriminant of this second order inequality is
\begin{equation}  \label{Delta1}
\Delta(l,h)=(19+8l\|\varphi\|)^{2}-32(3l\|\varphi\|+4l^{2}\|\varphi\|^{2}-%
\varepsilon).
\end{equation}
For $h,l$ small enough, this is estimated as
\begin{equation*}
\Delta(l,h)\sim361+32\varepsilon>0.
\end{equation*}%
Thus there are two zeros of the second order equality above $\eta_{1}=%
\displaystyle\frac{\sqrt{\Delta(l,h)}-(19+8l\|\varphi\|)}{16}>0$ and a
second zero $\eta_2<0$ rejected. Consequently, choosing $\eta\in]0,\eta_{1}[$
we obtain (\ref{LyapunovStability3-4}). Finally, (\ref{LyapunovStability3-3}%
) yields that $\|U^{1}\|\leq\varepsilon$. Next, equation (\ref%
{LyapunovStability4}), for $n=0$, implies that
\begin{equation}  \label{LyapunovStability4-1}
\|V^{1}\|\leq\,A(l,h,\varphi)\|U^{0}\|^{2}+B(l,h,\varphi)\|U^{0}\|+C(l,h,%
\varphi)+16|c_{2}|\|V^{0}\|,
\end{equation}%
where
\begin{equation*}
A(l,h,\varphi)=3+32|c_{2}|,
\end{equation*}%
\begin{equation*}
B(l,h,\varphi)=4\left(|c_{1}|+8|c_{2}|(2+l\|\varphi\|)+l\|\varphi\|+%
\displaystyle\frac{1}{h^{2}}\right) ,
\end{equation*}%
and
\begin{equation*}
C(l,h,\varphi)=2(1+8|c_{2}|)l^{2}\|\varphi\|^{2}+4l(4|c_{2}|+\displaystyle%
\frac{1}{h^{2}})\|\varphi\|.
\end{equation*}%
Choosing $\|(U^{0},V^{0})\|\leq\eta$, it suffices to study the inequality
\begin{equation}
A(l,h,\varphi)\eta^{2}+\left(B(l,h,\varphi)+16|c_{2}|\right)\eta+C(l,h,%
\varphi)-\varepsilon\leq0.  \label{V1Bound}
\end{equation}%
Its discriminant satisfies for $h,l$ small enough,
\begin{equation}  \label{Delta2}
\Delta(l,h)\sim\displaystyle\frac{16}{h^{4}}\left(1+20\alpha+|1-2\alpha|%
\right)^{2}+\displaystyle\frac{128\alpha|q|}{h^{2}}\varepsilon>0.
\end{equation}
Here also there are two zeros, $\eta_{1}^{\prime}=\displaystyle\frac{\sqrt{%
\Delta(l,h)}-(B(l,h,\varphi)+16|c_{2}|)}{2A(l,h,\varphi)}>0$ and a second
one $\eta^{\prime}<0$ and thus rejected. As a consequence, for $%
\eta\in]0,\eta _{1}^{\prime}[$ we obtain $\|V^{1}\|\leq\varepsilon$.
Finally, for $\eta\in]0,\eta_{0}[$ with $\eta_{0}=\min(\eta_{1},\eta_{1}^{%
\prime})$, we obtain $\|(U^{1},V^{1})\|\leq\varepsilon$ whenever $%
\|(U^{0},V^{0})\|\leq\eta$. Assume now that the $(U^{k},V^{k})$ is bounded
for $k=1,2,\dots,n$ (by $\varepsilon_{1}$) whenever $(U^{0},V^{0})$ is
bounded by $\eta$ and let $\varepsilon>0$. We shall prove that it is
possible to choose $\eta$ satisfying $\|(U^{n+1},V^{n+1})\|\leq\varepsilon$.
Indeed, from (\ref{LyapunovStability3-1}), we have
\begin{equation}  \label{LyapunovStability3Ordren-1}
\|U^{n+1}\|\leq19\varepsilon_{1}+8\varepsilon_{1}^{2}.
\end{equation}%
So, one seeks, $\varepsilon _{1}$ for which $8\varepsilon_{1}^{2}+19%
\varepsilon _{1}-\varepsilon \leq 0$. Its discriminant $\Delta=361+32%
\varepsilon $, with one positive zero $\varepsilon _{1}=\displaystyle\frac{%
\sqrt{361+32\varepsilon }-19}{16}$. Then $\|U^{n+1}\|\leq \varepsilon $
whenever $\|(U^{k},V^{k})\|\leq\varepsilon_{1}$, $k=1,2,\dots,n$. Next,
using (\ref{LyapunovStability4}) and (\ref{LyapunovStability3Ordren-1}), we
have
\begin{equation}
\|V^{n+1}\|\leq\left(4|c_{1}|+80|c_{2}|+1\right)\varepsilon_{1}+%
\left(32|c_{2}|+2\right)\varepsilon_{1}^{2}.
\label{LyapunovStability3Ordren-2}
\end{equation}%
So, it suffices as previously to choose $\varepsilon_{1}$ such that
\begin{equation*}
\left(32|c_{2}|+2\right)\varepsilon_{1}^{2}+\left(4|c_{1}|+80|c_{2}|+1%
\right)\varepsilon_{1}-\varepsilon\leq 0.
\end{equation*}%
$\Delta=(4|c_{1}|+80|c_{2}|+1)^{2}+4(32|c_{2}|+2)\varepsilon$, with positive
zero $\varepsilon_{1}^{\prime}=\displaystyle\frac{\sqrt{\Delta}%
-(4|c_{1}|+80|c_{2}|+1)}{2(32|c_{2}|+2)}$. Then $\|V^{n+1}\|\leq\varepsilon $
whenever $\|(U^{k},V^{k})\|\leq\varepsilon_{1}^{\prime}$, $k=1,2,\dots ,n$.
Next, it holds from the recurrence hypothesis for $\varepsilon _{0}=\min
(\varepsilon _{1},\varepsilon_{1}^{\prime })$, that there exists $\eta >0$
for which $\|(U^{0},V^{0})\|\leq \eta $ implies that $\|(U^{k},V^{k})\|\leq%
\varepsilon _{0}$, for $k=1,2,\dots ,n$, which by the next induces that $%
\|(U^{n+1},V^{n+1})\|\leq \varepsilon $.

\begin{lemma}
\label{laxequivresult} As the numerical scheme is consistent and stable, it
is then convergent.
\end{lemma}

This lemma is a consequence of the well known Lax-Richtmyer equivalence
theorem, which states that for consistent numerical approximations,
stability and convergence are equivalent. Recall here that we have already
proved in (\ref{consistency1}) and (\ref{consistency2}) that the used scheme
is consistent. Next, Lemma \ref{LyapunovStabilityLemma}, Lemma \ref%
{LWABounded} and equation (\ref{LyapunovStability1}) yields the stability of
the scheme. Consequently, the Lax equivalence Theorem guarantees the
convergence. So as Lemma \ref{laxequivresult}.

\section{Numerical implementations}

We propose in this section to present some numerical examples to validate
the theoretical results developed previously. The error between the exact
solutions and the numerical ones via an $L_{2}$ discrete norm will be
estimated. The matrix norm used will be
\begin{equation*}
\|X\|_{2}=\left( \displaystyle\sum_{i,j=1}^{N}|X_{ij}|^{2}\right)^{1/2}
\end{equation*}
for a matrix $X=(X_{ij})\in \mathcal{M}_{N+2}\mathbb{C}$. Denote $u^{n}$ the
net function $u(x,y,t^{n})$ and $U^{n}$ the numerical solution. We propose
to compute the discrete error
\begin{equation}
Er=\displaystyle\max_{n}\|U^{n}-u^{n}\|_{2}  \label{Er}
\end{equation}%
on the grid $(x_{i},y_{j})$, $0\leq\,i,j\leq\,J+1$ and the relative error
between the exact solution and the numerical one as
\begin{equation}
Relative\,Er=\displaystyle\max_{n}\displaystyle\frac{\|U^{n}-u^{n}\|_{2}}{%
\|u^{n}\|_{2}}  \label{Errelative}
\end{equation}%
on the same grid.

\subsection{A Polynomial-Exponential Example}

We develop in this part a classical example based on polynomial function
with an exponential envelop. We consider the inhomogeneous problem
\begin{equation}
\left\{
\begin{array}{l}
u_{tt}=\Delta\,u+v_{xx}+g(x,y,t),\quad(x,y,t)\in\Omega\times(t_{0},T), \\
v=qu_{xx}+u^{2},\quad(x,y,t)\in\Omega\times(t_{0},T), \\
\left(u,v\right)(x,y,t_{0})=\left(u_{0},v_{0}\right)(x,y),\quad(x,y,t)\in%
\overline{\Omega}\times(t_{0},T), \\
\frac{\partial\,u}{\partial\,t}(x,y,t_{0})=\varphi(x,y),\quad(x,y)\in%
\overline{\Omega}, \\
\overrightarrow{\nabla}(u,v)=0,\quad(x,y,t)\in\partial\Omega\times(t_{0},T)%
\end{array}
\right.
\end{equation}
where $\Omega=[-1,1]^{2}$ and where the right hand term is
\begin{eqnarray*}
g\left(x,y,t\right)&=&\left[\left(x^{2}-1\right)^{2}(x^{4}-58x^{2}+9)-48%
\left(35x^{4}-30x^{2}+3\right)+y^{4}-14y^{2}+5\right]e^{-t} \\
&&-16(x^{2}-1)^{2}\left[\left(x^{2}-1\right)^{4}\left(15x^{2}-1\right)+%
\left(y^{2}-1\right)^{2}\left(7x^{2}-1\right)\right]e^{-2t}
\end{eqnarray*}
The exact solution is
\begin{equation}
u(x,y,t)=\left[\left(x^{2}-1\right)^{4}+\left(y^{2}-1\right)^{2}\right]%
e^{-t}.
\end{equation}
In the following tables, numerical results are provided. We computed for
different space and time steps the discrete $L_{2}$-error estimates defined
by (\ref{Er}). The time interval is $[0,1]$ for a choice $t_{0}=0$ and $T=1$%
. The following results are obtained for different values of the parameters $%
J$ (and thus $h$), $l$ ((and thus $N$). The parameters $q$ and $\alpha $ are
fixed to $q=0.01$ and $\alpha =0.25$. We just notice that some variations
done on these latter parameters have induced an important variation in the
error estimates which explains the effect of the parameter $q$ which has the
role of a viscosity-type coefficient and the barycenter parameter $\alpha $
which calibrates the position of the approximated solution around the exact
one. Finally, some comparison with our work in \cite{Benmabrouk1} has proved
that Lyapunov type operators already result in fast convergent algorithms
with a maximum time of execution of 2.014 sd for the present one. The
classical tri-diagonal algorithms associated to the same problem with the
same discrete scheme and the same parameters yielded a maximum time of
552.012 sd, so a performance of $23.10^{-4}$ faster algorithm for the
present one. We recall that the tests are done on a Pentium Dual Core CPU
2.10 GHz processor and 250 Mo RAM.
\begin{table}[th]
\begin{center}
\centerline{Table 1.}
\begin{tabular}{||l|l|l|l||}
\hline\hline
J & $l$ & $Er$ & $Relative\,Er$ \\ \hline
10 & 1/100 & $4,0.10^{-3}$ & 0,1317 \\ \hline
16 & 1/120 & $3,3.10^{-3}$ & 0,1323 \\ \hline
20 & 1/200 & $2,0.10^{-3}$ & 0,1335 \\ \hline
24 & 1/220 & $1,8.10^{-3}$ & 0,1337 \\ \hline
30 & 1/280 & $1,4.10^{-3}$ & 0,1340 \\ \hline
40 & 1/400 & $9,8.10^{-4}$ & 0.1344 \\ \hline
50 & 1/500 & $7,8.10^{-4}$ & 0,1346 \\ \hline\hline
\end{tabular}%
\end{center}
\end{table}

\subsection{A 2-Particles Interaction Example}

The example developed hereafter is a model of interaction of two particles
or two waves. We consider the inhomogeneous problem
\begin{equation}
\left\{
\begin{array}{l}
u_{tt}=\Delta\,u+v_{xx}+g(x,y,t),\quad(x,y,t)\in\Omega\times(t_{0},T), \\
v=qu_{xx}+u^{2},\quad(x,y,t)\in\Omega\times(t_{0},T), \\
\left(u,v\right)(x,y,t_{0})=\left(u_{0},v_{0}\right)(x,y),\quad(x,y,t)\in%
\overline{\Omega}\times(t_{0},T), \\
\frac{\partial\,u}{\partial\,t}(x,y,t_{0})=\varphi(x,y),\quad(x,y)\in%
\overline{\Omega}, \\
\overrightarrow{\nabla}(u,v)=0,\quad(x,y,t)\in\partial\Omega\times(t_{0},T)%
\end{array}
\right.
\end{equation}
where
\begin{equation*}
g\left(x,y,t\right)=(4-6\psi^{2}(y))u^{2}-\psi^{2}(x)u.
\end{equation*}
and $u$ is the exact solution given by
\begin{equation*}
u(x,y,t)=2\psi^{2}(x)\psi^{2}(y)\theta(t)
\end{equation*}
with
\begin{equation*}
\psi(x)=\cos\left(\frac{x}{2}\right)\quad,\qquad\theta(t)=e^{-it}
\end{equation*}
and
\begin{equation*}
\varphi(x,y)=-2i\psi^{2}(x)\psi^{2}(y)
\end{equation*}
As for the previous example, the following tables shows the numerical
computations for different space and time steps the discrete $L_{2}$-error
estimates defined by (\ref{Er}) and the relative error (\ref{Errelative}).
The time interval is $[-2\pi,+2\pi]$ for a choice $t_{0}=0$ and $T=1$. The
following results are obtained for different values of the parameters $J$
(and thus $h$), $l$ ((and thus $N$). The parameters $q$ and $\alpha $ are
fixed here-also the same as previously, $q=0.01$ and $\alpha =0.25$.
Compared to the tri-diagonal scheme the present one leads a faster
convergent algorithms
\begin{table}[th]
\begin{center}
\centerline{Table 2.}
\begin{tabular}{||l|l|l|l||}
\hline\hline
J & $l$ & $Er$ & $Relative\,Er$ \\ \hline
10 & 1/100 & $4,6.10^{-3}$ & 0,2311 \\ \hline
16 & 1/120 & $4,4.10^{-3}$ & 0,2372 \\ \hline
20 & 1/200 & $2,4.10^{-3}$ & 0,2506 \\ \hline
24 & 1/220 & $2,3.10^{-3}$ & 0,2671 \\ \hline
30 & 1/280 & $2,0.10^{-3}$ & 0,3074 \\ \hline
40 & 1/400 & $1,4.10^{-3}$ & 0,3592 \\ \hline
50 & 1/500 & $7,6.10^{-4}$ & 0,2355 \\ \hline\hline
\end{tabular}%
\end{center}
\end{table}

\begin{remark}
For the convenience of the paper, we give here some computations of the
determinants $\Delta(l,h)$ for different values of the parameters of the
discrete scheme Firstly, for both examples above, we can easily see that $%
\|\varphi\|=$ and thus, equation (\ref{Delta1}) yields that
\begin{equation*}
\Delta(l,h)=361+32\varepsilon+416\,l-256\,l^2.
\end{equation*}
For the different values of $l$ as in the tables 1 and 2, we obtain a
positive discriminant leading two zeros with a rejected one. For the
discriminant of equation (\ref{Delta2}) we obtain
\begin{equation*}
\Delta(l,h)=\displaystyle\frac{676}{h^4}+\displaystyle\frac{8\varepsilon}{h^2%
}.
\end{equation*}
Hence, the results explained previously hold.
\end{remark}

\section{Conclusion}

This paper investigated the solution of the well-known Boussinesq equation
in two-dimensional case by applying a two-dimensional finite difference
discretization. The Boussinesq equation in its original form is a 4-th order
partial differential equation. Thus, in a first step it was recasted into a
system of second order partial differential equations using a reduction
order idea. Next, the system has been transformed into an algebraic discrete
system involving Lyapunov-Syslvester matrix terms by using a full time-space
discretization. Solvability, consistency, stability and convergence are then
established by applying well-known methods such as Lax-Richtmyer equivalence
theorem and Lyapunov Stability and by examining the Lyapunov-Sylvester
operators. The method was finally improved by developing numerical examples.
It was shown to be efficient by means of error estimates as well as time
execution algorithms compared to classical ones.

\section{Appendix}

\subsection{The Tridiagonal Associated System}

Consider the lexicographic mesh $k=j(J+1)+m$ for $0\leq\,j,m\leq\,J$, and
denote $N=J(J+2)$, and
\begin{equation*}
\Lambda_N=\{nJ+n-1\;;\;n\in\mathbb{N}\}\;\;,\;\;\widetilde{\Lambda}%
_N=\{n(J+1)\;;\;n\in\mathbb{N}\}\quad\hbox{and}\quad \Theta_N=\Lambda_N\cup%
\widetilde{\Lambda}_N.
\end{equation*}
We obtain a tri-diagonal block system on the form
\begin{equation}  \label{tridiagonalsystem1}
\left\{%
\begin{matrix}
\medskip\widetilde{A}U^{n+1}+a_2\widetilde{R}V^{n+1}=\widetilde{B}U^{n}-%
\widetilde{A}U^{n-1}+b_2\widetilde{R}V^{n}-a_2\widetilde{R}V^{n-1}\hfill\cr%
\medskip \medskip\,V^{n+1}-2c_2\widetilde{R}U^{n+1}=2c_1\widetilde{R}%
U^{n}+2c_2\widetilde{R}U^{n-1}-V^{n-1}+2\widehat{F^{n}}.\hfill%
\end{matrix}%
\right.
\end{equation}
The numerical solutions' matrices $U^n$ and $V^n$ are identified here as
one-column $(N+1)$-vectors and the matrices $\widetilde{A}$, $\widetilde{B}$
and $\widetilde{R}$ are evaluated as follows.\newline
\underline{The matrix $\widetilde{A}$}
\begin{equation*}
\widetilde{A}_{j,j}=2a_1\;,\;\forall\,j\;\;\;;0\leq\,j\leq\,N,
\end{equation*}
\begin{equation*}
\widetilde{A}_{j,j+1}=\displaystyle\frac{1}{2}\widetilde{A}%
_{0,1}=a_2\;,\;\forall\,j\;;\;\;1\leq\,j\leq\,N\,;\;\;
j\notin\Theta_{N},\quad\hbox{and}\quad0\;\;\hbox{on}\;\;\Lambda_N,
\end{equation*}
\begin{equation*}
\widetilde{A}_{j-1,j}=\displaystyle\frac{1}{2}\widetilde{A}%
_{N,N-1}=a_2\;,\;\forall\,j\;;\;\;1\leq\,j\leq\,N\,;\;\;
j\notin\Theta_{N},,\quad\hbox{and}\quad0\;\;\hbox{on}\;\;\widetilde{\Lambda}%
_N,
\end{equation*}
\begin{equation*}
\widetilde{A}_{j,j+J+1}=2a_2\;,\;\forall\,j\;,\;\;0\leq\,j\leq\,J
\end{equation*}
\begin{equation*}
\widetilde{A}_{j,j+J+1}=a_2\;,\;\forall\,j\;,\;\;J+1\leq\,j\leq\,N-J-1
\end{equation*}
\begin{equation*}
\widetilde{A}_{j-J-1,j}=a_2\;,\;\forall\,j\;,\;\;J+1\leq\,j\leq\,N-J-1,
\end{equation*}
\begin{equation*}
\widetilde{A}_{j-J-1,j}=2a_2\;,\;\forall\,j\;,\;\;N-J\leq\,j\leq\,N.
\end{equation*}
\underline{The matrix $\widetilde{B}$}
\begin{equation*}
\widetilde{B}_{j,j}=2b_1\;,\;\forall\,j\;\;\;;0\leq\,j\leq\,N,
\end{equation*}
\begin{equation*}
\widetilde{B}_{j,j+1}=\displaystyle\frac{1}{2}\widetilde{B}%
_{0,1}=b_2\;,\;\forall\,j\;;\;\;1\leq\,j\leq\,N\,;\;\;
j\notin\Theta_{N},\quad\hbox{and}\quad0\;\;\hbox{on}\;\;\Lambda_N,
\end{equation*}
\begin{equation*}
\widetilde{B}_{j-1,j}=\displaystyle\frac{1}{2}\widetilde{B}%
_{N,N-1}=b_2\;,\;\forall\,j\;;\;\;1\leq\,j\leq\,N\,;\;\;
j\notin\Theta_{N},,\quad\hbox{and}\quad0\;\;\hbox{on}\;\;\widetilde{\Lambda}%
_N,
\end{equation*}
\begin{equation*}
\widetilde{B}_{j,j+J+1}=2b_2\;,\;\forall\,j\;,\;\;0\leq\,j\leq\,J
\end{equation*}
\begin{equation*}
\widetilde{B}_{j,j+J+1}=b_2\;,\;\forall\,j\;,\;\;J+1\leq\,j\leq\,N-J-1
\end{equation*}
\begin{equation*}
\widetilde{B}_{j-J-1,j}=b_2\;,\;\forall\,j\;,\;\;J+1\leq\,j\leq\,N-J-1,
\end{equation*}
\begin{equation*}
\widetilde{B}_{j-J-1,j}=2b_2\;,\;\forall\,j\;,\;\;N-J\leq\,j\leq\,N.
\end{equation*}
\underline{The matrix $\widetilde{R}$}
\begin{equation*}
\widetilde{R}_{j,j}=-2\;,\;\forall\,j\;\;\;;0\leq\,j\leq\,N,
\end{equation*}
\begin{equation*}
\widetilde{R}_{j,j+J+1}=2\;,\;\forall\,j\;,\;\;0\leq\,j\leq\,J.
\end{equation*}
\begin{equation*}
\widetilde{R}_{j,j-J-1}=2\;,\;\forall\,j\;,\;\;N-J\leq\,j\leq\,N.
\end{equation*}
\begin{equation*}
\widetilde{R}_{j,j+J+1}=\widetilde{R}_{j-J-1,j}=1\;,\;\forall\,j\;,\;\;J+1%
\leq\,j\leq\,N-J-1.
\end{equation*}

\subsection{Headlines of the algorithm applied}

\begin{itemize}
\item Compute the matrices of the system

\item Initialisation: Compute the matrices $U^0$, $U^1$, $V^0$ and $V^1$

\item for $n\geq2$,
\begin{equation*}
U^n=lyap(W,A,\mathcal{L}_{\widetilde{B},B}(U^{n-1})+b_2RV^{n-1}-\mathcal{L}%
_{W,A}(U^{n-2})-a_2R(F^{n-2}+F^{n-1})),
\end{equation*}
and
\begin{equation*}
V^{n+1}=2c_2RU^{n}+2R(c_1U^{n-1}+c_2U^{n-2})-V^{n-2}+2\widehat{F^{n-1}}.
\end{equation*}
\end{itemize}


\begin{thebibliography}{99}
\bibitem{BaoDan1} T. Bao-Dan, Q. Yan-Hong and Chen-Ning, Exact Solutions for
a Class of Boussinesq Equation, Applied Mathematical Sciences, 3(6) (2009),
257-265.

\bibitem{Benmabrouk2} A. Ben Mabrouk and M. Ayadi, A linearized
finite-difference method for the solution of some mixed concave and convex
non-linear problems, Appl. Math. Comput. 197 (2008), 1--10.

\bibitem{Benmabrouk1} A. Ben Mabrouk and M. Ayadi, Lyapunov type operators
for numerical solutions of PDEs, Appl. Math. Comput. 204 (2008), 395--407.

\bibitem{Benmabrouk3} A. Ben Mabrouk, M.L. Ben Mohamed, K. Omrani,
Finite-difference approximate solutions for a mixed sub-superlinear
equation, Appl. Math. Comput. 187 (2007) 1007\^{a}\euro ``1016.

\bibitem{Bratsos1} A. G. Bratsos, Ch. Tsituras and D. G. Natsis, Linearized
numerical schemes for the Boussinesq equation. Appl. Num. Anal. Comp. Math.
2(1) (2005), 34-53.

\bibitem{Clarkson1} P. A. Clarkson, Nonclassical symmetry reductions for the
Boussinesq equation, Chaos, Solitons and Fractals, 5 (1995), 2261--2301.

\bibitem{Dehghan2006} M. Dehghan, Finite difference procedures for solving a
problem arising in modeling and design of certain optoelectronic devices,
Mathematics and Computers in Simulation 71 (2006), 16--30.

\bibitem{Dehghan2008} M. Dehghan, Use of Hes homotopy
perturbation method for solving a partial differential equation arising in
modeling of flow in porous media, Journal of Porous Media, Volume 11(8)
(2008), 765--778.

\bibitem{Dehghan-Salehi} M. Dehghan, and R. Salehi, A meshless based
numerical technique for traveling solitary wave solution of Boussinesq
equation, Applied Mathematical Modelling 36 (2012), 1939--1956.

\bibitem{ElMikkawy1} M. El-Mikkawy, A note on a three-term recurrence for a
tridiagonal matrix, Appl. Math. Computa., 139 (2003), 503-511.

\bibitem{ElMikkawy2} M. El-Mikkawy, A fast algorithm for evaluating $n$th
order tri-diagonal determinants, J. Computa. \& Appl. Math., 166 (2004),
581-584.

\bibitem{ElMikkawy3} M. El-Mikkawy, On the inverse of a general tridiagonal
matrix, J. Computa. \& Appl. Math., 150 (2004), 669-679.

\bibitem{ElMikkawy-Karawia} M. El-Mikkawy and A. Karawia, Inversion of
general tridiagonal matrices, Appl. Math. Letters 19 (2006), 712-720.

\bibitem{El-Mikkawy-Atlan} M. El-Mikkawy, F. Atlan, A new recursive
algorithm for inverting general image-tridiagonal matrices. Applied
Mathematics Letters, 44 (2015), 34-39.

\bibitem{Goncalves} E. Gon\c calv\`es, Resolution numerique, discretisation
des EDP et EDO. Cours, Institut National Polytechnique de Grenoble, 2005.

\bibitem{HenriCartan} H. Cartan, Differential Calculus, Kershaw Publishing
Company, London 1971, Translated from the original French text Calcul
differentiel, first published by Hermann in 1967.

\bibitem{Jafarizadeh1} M. A. Jafarizadeh, A. R. Esfandyari and M.
Moslehi-fard, Exact solutions of Boussinesq equation. Third International
Conference on Geometry, Integrability and Quantization June 14.23, 2001,
Varna, Bulgaria Ivailo M. Mladenov and Gregory L. Naber, Editors Coral
Press, Sofia 2001, pp 304-314.

\bibitem{Jameson} A. Jameson, Solution of equation $AX+XB=C$ by inversion of
an $M\times M$ or $N\times N$ matrix. SIAM J. Appl Math. 16(5) (1968),
1020-1023.

\bibitem{Jia-Li} J. Jia, S. Li, On the inverse and determinant of general
bordered tridiagonal matrices. Computers \& Mathematics with Applications,
69(6) (2015), 503-509.

\bibitem{Kano1} T. Kano and T. Nishida, A mathematical justification for
Korteweg-De Vries equation and Boussinesq equation of water surface waves.
Osaka J. Math. 23 (1986), 389-413.

\bibitem{Kaya1} D. Kaya, Explicit solutions of generalized nonlinear
Boussinesq equations. Journal of Applied Mathematics 1(1) (2001), 29--37.

\bibitem{Lai1} Sh. Lai, Y. Wub and Y. Zhou, Some physical structures for the
(2+1)-dimensional Boussinesq water equation with positive and negative
exponents. Computers and Mathematics with Applications 56 (2008), 339--345.

\bibitem{Liu1} Ch. Liu and Z. Dai, Exact periodic solitary wave solutions
for the (2+1)-dimensional Boussinesq equation. J. Math. Anal. Appl. 367
(2010), 444--450.

\bibitem{Parlange1} J.-Y. Parlange, W. L. Hogarth, R. S. Govindaraju, M. B.
Parlange and D. Lockington, On an Exact Analytical Solution of the
Boussinesq Equation, Transport in Porous Media 39 (2000), 339--345.

\bibitem{Rionero} S. Rionero, Stability results for hyperbolic and parabolic
equations. Transport Theory \& Statistical Physics, 25(3-5) (1996), 323-337.

\bibitem{Shokri-Dehghan} A. Shokri, and M. Dehghan, A Not-a-Knot meshless
method using radial basis functions and predictor-corrector scheme to the
numerical solution of improved Boussinesq equation, Computer Physics
Communications 181 (2010), 1990--2000.

\bibitem{Song1} M. Song and Sh. Shao, Exact solitary wave solutions of the
generalized (2+1)-dimensional Boussinesq equation. Applied Mathematics and
Computation 217 (2010), 3557--3563.

\bibitem{Varlamov1} V. Varlamov, Two-dimensional Boussinesq equation in a
disc and anisotropic Sobolev spaces. C. R. Mecanique 335 (2007), 548--558.

\bibitem{Wazwaz1} A.-M. Wazwaz, Variants of the Two-Dimensional Boussinesq
Equation with Compactons, Solitons, and Periodic Solutions. Computers and
Mathematics with Applications 49 (2005), 295-301.

\bibitem{Yi1} Z. Yi, Y. Ling-ya, L. Yi-neng and Z. Hai-qiong, Periodic wave
solutions of the Boussinesq equation. J. Phys. A: Math. Theor. 40(21)
(2007), 5539--5549.

\bibitem{Yi2} Z. Yi and Y. Ling-Ya, Rational and Periodic Wave Solutions of
Two-Dimensional Boussinesq Equation. Commun. Theor. Phys. 49(4) (2008),
815-824.
\end{thebibliography}
\end{document}